\documentclass[12pt]{article}
\usepackage{amsmath,amssymb,amsthm,latexsym}

\newtheorem{Def}{Definition}[section]
  \newtheorem{Thm}[Def]{Theorem}
  \newtheorem{Lem}[Def]{Lemma}
  \newtheorem{Prop}[Def]{Proposition}
  
  \newtheorem{Rem}[Def]{Remark}
\newcommand{\Hom}{\mathop{\mathrm{Hom}}\nolimits}

\numberwithin{equation}{section}

\begin{document}
\begin{center}
\large{\textbf{Symmetric linear functions
of
the restricted quantum group $\overline{U}_qsl_2(\mathbb{C})$}}
\end{center}

\begin{center}
Yusuke Arike \\
\small{Department of Pure and Applied Mathematics, \\ 
Graduate School of Imformation Science and Technology, \\
Osaka University,
Toyonaka, Osaka 560-0043, Japan \\
E-mail: y-arike@cr.math.sci.osaka-u.ac.jp }
\end{center}

\abstract
We determine a set of primitive idempotents and the basic algebra of the restricted quantum group $\overline{U}_qsl_2(\mathbb{C})$.
As a result, we can show the dimension of the space of symmetric linear functions of $\overline{U}_qsl_2(\mathbb{C})$
is $3p-1$.

\section{Introduction}
In the study of conformal field theory associated with a vertex operator algebra(VOA),
the representation theory of VOAs plays an important role.
In fact the theory for any rational vertex operator algebra with  the factorization property
is established over the projective line in \cite{NT}. This theory is generalized to the higher genus case;
in particular, for an elliptic curve the space of conformal blocks with a vacuum module $V$
is nothing but the space of formal characters of modules for $V$ (cf. \cite{Z}).
Therefore its dimension coincides with the number of simple modules for $V$ up to
isomorphisms.

On the other hand, the theory for irrational VOAs is rather difficult.  For exapmle,
under the same finiteness conditions it is shown that genus one conformal blocks
is finite dimensional (see \cite{Mi} ). In this case the dimension is greater than 
the number of simple modules. 

There is slightly well studied example of irrational conformal field theory
which is called logarithmic conformal field theory. A typical example of theory
is a $\mathcal{W}(p)$-algebra, whose confomal blocks involve logarithmic function 
of  modulus $q$ (compare that  no logarithmic terms appear in rational cases.) 
In these examples determining the dimension of conformal blocks is
not easy, in fact,  this dimension  has not been yet known.

The discussion given in \cite{Mi} and a private communication \cite{Na}
suggest that this dimension coincides with the one of the symmetric linear functions
of a finite dimensional algebra whose category of modules is equivalent to the category
of $V$-modules (also see \cite{MNT}).

Recently \cite{FGST1} and \cite{FGST2} found that the category of $\mathcal{W}(p)$-modules
is closely related with the category of the restricted quantum group
$\overline{U}_qsl_2(\mathbb{C})$(for short $\overline{U}_q$) with $q=e^{\pi i/p}$.
More precisely they proved that
a subspace of conformal blocks of $W(p)$-algebras has a subspace
which is invariant under the canonical $SL_2(\mathbb{Z})$ action
and this space is isomorphic to the space of $q$-characters.
We hope that the space of $q$-characters is isomorphic to conformal blocks.
They further proved that if $p=2$ the category of modules for $\overline{U}_q$
and the category of modules for $\mathcal{W}(p)$. We can expect that this equivalence of 
categories holds for all $p \ge 2$.

The above situation naturally indicates that the dimension of conformal blocks is exactly the dimension of $\overline{U}_q/[\overline{U}_q,\overline{U}_q]$, that is, the dimension of the space of the symmetric linear functions over the algebra
$\overline{U}_q$.
This is still conjecture, but assuming 
that this is true, we are going to determine the dimension of $\overline{U}_q/[\overline{U}_q,\overline{U}_q]$.

It is proved in \cite{FGST1} that $\overline{U}_q = \bigoplus_{s=0}^{p} \mathcal{Q}_s $ where $\mathcal{Q}_s$
is a subalgebra of $\overline{U}_q$.
Thus it suffices to determine the dimension of the space $\mathcal{Q}_s/[\mathcal{Q}_s, \mathcal{Q}_s]$.
Nesbitt and Scott showed in \cite{NS} that the space of symmetric linear functions of a finite-dimensional algebra is isomorphic to the one of its basic algebra as a vector space.
We prove that the basic algebra $B_s$ of $\mathcal{Q}_s$ for $1 \le s \le p-1$ is $8$-dimensional and that $B_0$ and $B_p$ are $1$-dimensional.
Moreover we get the multiplication table among the basis of $B_s$ for $1 \le s \le p-1$ and this implies that $\dim B_s/[B_s,B_s] = 3$. 
Therefore the space of symmetric linear functions of $\overline{U}_q$ is $(3p-1)$-dimensional. 

This paper is organized as follows.
In section $2$ we introduce the basic algebra of a finite-dimensional associative algebra,
symmetric linear functions and the result given by Nesbitt and Scott
which states that the dimension of the space of symmetric linear functions of an algebra is equal to 
that of its basic algebra.
In section $3$ we recall some definitions and results in the restricted quantum group $\overline{U}_q$
studied(\cite{FGST1}, \cite{FGST2}).
Moreover we determine a basic set of primitive idempotents $\{ e_s^{\pm} | 1 \le s \le p  \}$of $\overline{U}_q$.
The primitive idempotent $e_s^{\pm}$ generates the indecomposable projective module $\mathcal{P}_s^{\pm}$
and $\mathcal{P}_s^{\pm}$ is the $2p$-dimensional module whose
basis is given by $\{ B_n^{\pm}(s), X_k^{\pm}(s), Y_k^{\pm}(s), A_n^{\pm}(s) | 0 \le n \le s-1, 0 \le k \le p-s-1 \}$
where $B_0^{\pm}(s) = e_s^{\pm}$.
In section $4$ we prove that the basic algebra of the subalgebra $\mathcal{Q}_s$ for $1 \le s \le p-1$
of $\overline{U}_q$ is $8$-dimensional and $Q_0$ and $Q_p$ are $1$-dimensional
where $\overline{U}_q = \bigoplus_{s=0}^{p} \mathcal{Q}_s $.
And we find that the basis of $\mathcal{Q}_s$ for $1 \le s \le p-1$ is given by
$\{e_s^+, e_{p-s}^-, X_0^{\pm}, Y_0^{\pm}, A_0^{\pm} \}$ 
where $Z_0^+ = Z_0^+(s)$ and $Z_0^- = Z_0^-(p-s)$ for $Z = X, Y, A$.
We also determine the multiplication talbe among the basis of  $\mathcal{Q}_s$ for $1 \le s \le p-1$.
We prove that the space of symmetric linear functions of the baisc algebra of $\mathcal{Q}_s$ for $1 \le s \le p-1$
is $3$-dimensional.
Then we prove that the space of symmetric linear functions of $\overline{U}_q$ is $3p-1$-dimensional.

The author thanks to Professor K. Nagatomo for continuous encouragement.

\section{Preliminaries}
In this section, we introduce the basic algebra of a finite-dimensional associative algebra and
symmetric linear functions.
\subsection{The basic algebra}
Let $A$ be a finite-dimensional associative algebra over $\mathbb{C}$ with a unity $1$
and 
\begin{align}
1 = \sum_{i=1}^{n}\sum_{j=1}^{n_i} e_{ij}, \ Ae_{ij} \cong Ae_{ik}, \ Ae_{ij} \not\cong Ae_{kl} \ (i \neq k) \notag
\end{align}
be a decomposition into primitive idempotents of the unity.
We now set $e_i = e_{ij_a}$ for some $1 \le j_a \le n_i$ and call $\mathcal{E} = \{ e_i \ | \ 1 \le i \le n \}$
\textit{the basic set of primitive idempotents}.
Let $\mathcal{S}$ be the set of equivalence classes of all simple left $A$-modules
and $\mathcal{P}$ be the set of equivalence classes of all indecomposable projective left $A$-modules.
It is well known that each of the maps $\mathcal{E} \to \mathcal{S}$ defined by $e \mapsto Ae/Je$
and $\mathcal{E} \to \mathcal{P}$ definedby $e \mapsto Ae$ is bijective where $J$ is the Jacobson radical of $A$.

Set $e = \sum_{i=1}^{n} e_i$.
Then the subspace $B(A) = eAe$ is the associative algebra with the unity $e$,
which is called the \textit{basic algebra} of $A$.

\subsection{Symmetric linear functions}
Let $A$ be a finite-dimensional associative algebra with a unity $1$.
A \textit{symmetric linear function} of $A$ is a linear function $\varphi \in \Hom_{\mathbb{C}}(A, \mathbb{C})$
satisfying the condition $\varphi (ab) = \varphi (ba)$ for all $a, b \in A$.
This implies that the dimension of the space of the symmetric linear functions of $A$ is given by $\dim A/[A, A]$.

The following theorem can be found in \cite{NS}.

\begin{Thm}[\cite{NS}]\label{Thm:NS}
The space of symmetric linear functions of a fnite-dimensional algebra is isomorphic
to the one of its basic algebra.
\end{Thm} 

\section{The restricted quantum group $\overline{U}_qsl_2(\mathbb{C})$ and its modules}
In this section, we review the restricted quantum group $\overline{U}_q$
and its representations(\cite{FGST1}, \cite{FGST2}).
The most important result here is determining the basic set of primitive idempotents of $\overline{U}_q$.
We set $q = e^{\pi i /p}$ for an integer $p \ge 2$
and 
\begin{align}
[n] = \frac{q^n-q^{-n}}{q-q^{-1}} \notag
\end{align}
for any integer $n$.

\subsection{The restricted quantum group $\overline{U}_q$}
The restricted quantum group $\overline{U}_qsl_2(\mathbb{C})(= : \overline{U}_q)$ is an associative algebra generated by $E$, $F$ and $K$ with the following relations:
\begin{align}
E^p = F^p = 0,\ \ & K^{2p} = 1, \notag\\
KEK^{-1} = q^2E,\ \ KFK^{-1} =& q^{-2}F,\ \ [E,F] = \frac{K - K^{-1}}{q - q^{-1}}. \notag
\end{align}
Moreover this algebra has a structure of Hopf algebra,
though we do not use it in this paper(see \cite{FGST1}, \cite{FGST2}).

The classification of simple left $\overline{U}_q$-modules can be found in \cite{FGST1}.
\begin{Thm}[\cite{FGST1}]
The complete list of inequivalent simple modules is given as follows:
$\mathcal{X}_s^{\pm}$ with basis $\{a_n^{\pm}\}_{n=0}^{s-1}$ where $1 \le s \le p$.
the vector $a_0^{\pm}$ is a highest weight vector with highest weight $\pm q^{s-1}$, and the algebra action given by
\begin{align}
K a_n^{\pm} =& \pm q^{s-1-2n}a_n^{\pm}, \notag\\
E a_n^{\pm} =& \pm [n] [s-n]  a_{n-1}^{\pm},\notag\\
F a_n^{\pm} =&   a_{n+1}^{\pm}, \notag
\end{align}
where we set $a_{-1}^{\pm} = a_{s}^{\pm} = 0$.
\end{Thm}

\subsection{Indecomposable projective modules and the basic set of primitive idempotents}
In this section, we construct primitive idempotents of $\overline{U}_q$,
and classify the inequivalent indecomposable projective modules.
The following lemma will be often used in this paper.
\begin{Lem}[\cite{K}] \label{Lem:2}
For $1 \leq m \leq p - 1$, The following relations hold in $\overline{U}_q$:
\begin{align}
[E,F^m]  & =   [m]F^{m-1}\frac{q^{-(m - 1)}K - q^{m - 1}K^{-1}}{q - q^{-1}} \notag\\
         & =   [m]\frac{q^{m - 1}K - q^{- (m - 1)}K^{-1}}{q - q^{-1}}F^{m - 1} \notag\\
[E^m,F]  & =   [m]E^{m-1}\frac{q^{m - 1}K - q^{-(m - 1)}K^{-1}}{q - q^{-1}} \notag\\
         & =   [m]\frac{q^{-(m - 1)}K - q^{ m - 1}K^{-1}}{q - q^{-1}}E^{m - 1} \notag
\end{align}
\end{Lem}
Set $v_s^{\pm} = \sum_{\ell = 0}^{2p-1} (\pm q^{-(s-1)})^{\ell} K^{\ell} \in \overline{U}_q$ for $1 \le s \le p$.
Since $Kv_s^{\pm} = \pm q^{s-1} v_s^{\pm}$,
we see
$a_0^{\pm}(s) = E^{p-1}F^{p-1} v_s^{\pm}$ is the highest weight vector of highest weight $\pm q^{s-1}$.
We now construct simple modules as submodules of $\overline{U}_q$ (see Remark \ref{rem:1})
\begin{Lem}\label{proj:1}
Set $a_n^{\pm}(s) = F^n a_0^{\pm}(s)$,
then
\begin{align}
a_n^{\pm}(s) = \begin{cases}
               \prod_{i=1}^{n} (\pm [i] [s-i]) E^{p-1-n} F^{p-1}  v_s^{\pm}, & 1 \le n \le s-1, \\
               0, & n \ge s.
               \end{cases} \notag
\end{align}
\end{Lem}
\begin{proof}
By Lemma \ref{Lem:2}, we have
\begin{align}
FE^{p-1} = E^{p-1}F - \frac{[p-1]}{q-q^{-1}}E^{p-2}(q^{p-2}K - q^{-(p-2)}K^{-1}). \notag
\end{align}
This implies that
\begin{align}
 a_1^{\pm}(s) &= \left\{ E^{p-1}F - \frac{[p-1]}{q-q^{-1}}E^{p-2}(q^{p-2}K - q^{-(p-2)}K^{-1}) \right\} F^{p-1}  v_s^{\pm} \notag\\
&=  \frac{[1]}{q-q^{-1}}(\pm q^{s-1} \mp q^{-(s-1)} )E^{p-2}F^{p-1}  v_s^{\pm} \notag\\
&= \pm [1] [s-1] E^{p-2}F^{p-1}  v_s^{\pm}. \notag
\end{align}

For $n \ge 2$
Lemma \ref{Lem:2} shows the following:
\begin{align}
& F \prod_{i=1}^{n-1} (\pm [i] [s-i]) E^{p-n} F^{p-1}  v_s^{\pm} \notag\\
& =\prod_{i=1}^{n-1} (\pm [i] [s-i]) \left\{ E^{p-n}F - \frac{[p-n]}{q - q^{-1}}E^{p-1-n} (q^{p-1-n}K - q^{-(p-1-n)}K^{-1}) \right\} F^{p-1}v_s^{\pm} \notag\\
& = \prod_{i=1}^{n-1} (\pm [i] [s-i]) \pm [n][s-n] E^{p-1-n}F^{p-1} v_s^{\pm} \notag
\end{align}
Moreover we have
\begin{align}
&F E^{p-s} F^{p-1} v_s^{\pm} \notag\\
&= \left\{ E^{p-s}F - \frac{[p-s]}{q - q^{-1}}E^{p-s-1}(q^{p-s-1}K - q^{-(p-s-1)}K^{-1}) \right\}F^{p-1} v_s^{\pm} \notag\\
& = - \frac{[p-s]}{q - q^{-1}}E^{p-s-1}(\pm q^{p-s-1}  q^{s-1-2(p-1)} \mp q^{-(p-s-1)}  q^{-(s-1)+2(p-1)} ) F^{p-1} v_s^{\pm} \notag\\
&= -  \frac{[p-s]}{q - q^{-1}} (\pm q^{p} \mp q^{-p})E^{p-s-1}F^{p-1} v_s^{\pm} = 0. \notag
\end{align}
\end{proof}
\begin{Rem} \label{rem:1}
By the above lemma, we see 
\begin{align}
K a_n^{\pm}(s) =& \pm q^{s-1-2n}a_n^{\pm}(s), \notag\\
E a_n^{\pm}(s) =& \pm [n] [s-n]  a_{n-1}^{\pm}(s),\notag\\
F a_n^{\pm}(s) =&   a_{n+1}^{\pm}(s), \notag
\end{align}
and $E a_0^{\pm}(s) = 0$ and $F a_{s-1}^{\pm}(s)=0$.
So the subspace spanned by $\{a_n^{\pm}(s) | \ 0 \le n \le s-1 \}$
is isomorphic to the simple module $\mathcal{X}_s^{\pm}$.
\end{Rem}
We next determine a primitive idempotent.
\begin{Prop}\label{proj:2}
The element 
\begin{align}
e_p^{\pm} = \frac{1}{2p \prod_{i=1}^{p-1}(\pm[i][s-i])}a_0^{\pm}(p) \notag
\end{align} 
is a primitive idempotent of $\overline{U}_q$.
In paticular, the simple module $\mathcal{X}^{\pm}_{p}$ is projective.
\end{Prop}
\begin{proof}
This follows from
\begin{align}
(a_0^{\pm}(p) )^{2} &= E^{p-1}F^{p-1}v_p^{\pm} a_0^{\pm}(p) \notag\\
& = E^{p-1}F^{p-1} \sum_{\ell = 0}^{2p-1} (\pm q^{-(p-1)})^{\ell} (\pm q^{p-1})^{\ell} a_0^{\pm}(p) \notag\\
& = 2p\prod_{i = 1}^{p-1} (\pm [i][s-i])a_0^{\pm}(p). \notag
\end{align}
\end{proof}

Next we consider the case $1 \le s \le p-1$.
\begin{Lem}\label{proj:3}
For $1 \le s \le p-1$,
\begin{align}
a_0^{\pm}(s) = F\sum_{n=1}^{p-s} \alpha_n^{\pm}(s)E^{p-n}F^{p-1-n}v_s^{\pm} \notag
\end{align}
where
\begin{align}
\alpha_n^{\pm}(s)  = \prod_{k=p-s-(n-1)}^{p-s-1} (\mp [k][p-s-k]). \notag
\end{align}
\end{Lem}
\begin{proof}
By Lemma \ref{Lem:2},
we have
\begin{align}
&E^{p-n}F^{p-n}v_{s}^{\pm} \notag\\
&= \left\{ FE^{p-n} + \frac{[p-n]}{q - q^{-1}}E^{p-1-n}(q^{p-1-n}K - q^{-(p-1-n)}K^{-1}) \right\} F^{p-1-n}v_{s}^{\pm} \notag\\
&= FE^{p-n} F^{p-1-n}v_{s}^{\pm} \mp [n][p-s-n]E^{p-1-n}F^{p-1-n}v_{s}^{\pm} \notag
\end{align}
and $E^{s}F^{s}v_{s}^{\pm} = F E^{s}F^{s-1}v_{s}^{\pm}$.
Moreover the relation
\begin{align}
\alpha_n^{\pm}(s) = \mp[n-1][p-s-(n-1)]\alpha_{n-1}^{\pm}(s) \notag
\end{align}
shows 
\begin{align}
\alpha_n^{\mp}(s) = \prod_{k=1}^{n-1}(\mp[k][p-s-k]) = \prod_{k=p-s-(n-1)}^{p-s-1}(\mp[k][p-s-k]).   \notag 
\end{align}\qedhere 
\end{proof}

Now we set
\begin{align}
&b_0^{\pm}(s) = \sum_{n=1}^{p-s} \alpha_n^{\pm}(s)E^{p-1-n}F^{p-1-n}v_s^{\pm},\notag\\
&b_n^{\pm}(s) = F^nb_0^{\pm}(s),  \notag\\
&x_{k}^{\pm}(s)  = \frac{1}{\prod_{\ell = k+1}^{p-s-1}(\mp[\ell][p-s-\ell])}E^{p-s-k}b_{0}^{\pm}(s), \label{proj:re0}\\
&y_{k}^{\pm}(s)  = F^{s+k}b_0^{\pm}(s), \label{proj:re0-1}
\end{align}
for $1 \le n \le s-1$ and $1 \le k \le p-s-1$.
Let $\mathcal{P}_s^{\pm}$ be the subspace spanned by
\begin{align}
\{ b_n^{\pm}(s), a_n^{\pm}(s), x_{k}^{\pm}(s), y_{k}^{\pm}(s) \ | \ 0 \le n \le s-1, \ 0 \le k \le p-s-1 \}.\notag
\end{align}

\begin{Prop}\label{proj:4}
\begin{enumerate}
\item[$(1)$]
The subspace $\mathcal{P}_s^{\pm}$ of
$\overline{U}_q$ is the $2p$-dimensional left $\overline{U}_q$-module
with the basis
\begin{align}
b_n^{\pm}(s), \ a_n^{\pm}(s), \ x_{k}^{\pm}(s), \ y_{k}^{\pm}(s) \notag
\end{align}
for $0 \le n \le s-1$ and $0 \le k \le p-s-1$
and the action of $\overline{U}_q$ is given by
\begin{align}
&Kb_n^{\pm}(s) = \pm q^{s-1-2n}b_n^{\pm}(s), \ Ka_n^{\pm}(s) = \pm q^{s-1-2n}a_n^{\pm}(s),  \notag\\
&Kx_{k}^{\pm}(s) = \mp q^{p-s-1-2k}x_{k}^{\pm}(s), \ Ky_{k}^{\pm}(s) = \mp q^{p-s-1-2k}y_{k}^{\pm}(s), \notag
\end{align}
\begin{align}
&Eb_n^{\pm}(s) = \begin{cases}
                    \pm [n][s-n] b_{n-1}^{\pm}(s) + a_{n-1}^{\pm}(s), & 1 \le n \le s-1, \\
                    x_{p-s-1}^{\pm}(s), & n=0,
                    \end{cases}\notag\\
&Ex_{k}^{\pm}(s) = \begin{cases}
                    \mp [k][p-s-k]x_{k-1}^{\pm}(s), & 1 \le k \le p-s-1, \\
                    0, & k=0,
                    \end{cases}\notag
                    \end{align}
                    \begin{align}
&Ey_{k}^{\pm}(s) = \begin{cases}
                    \mp [k][p-s-k]y_{k-1}^{\pm}(s), & 1 \le k \le p-s-1, \\
                    a_{s-1}^{\pm}(s), & k=0,
                    \end{cases}\notag\\
&Ea_{n}^{\pm}(s) = \begin{cases}
                    \pm [n][s-n]a_{n-1}^{\pm}(s), & 1 \le n \le s-1, \\
                    0, & n=0,
                    \end{cases}\notag
                    \end{align}
                    \begin{align}
&Fb_n^{\pm}(s) = \begin{cases}
                 b_{n+1}^{\pm}(s), & 0 \le n \le s-2, \\
                 y_0^{\pm}(s), & n = s-1,
                 \end{cases}\notag\\
&Fx_{k}^{\pm}(s) = \begin{cases}
                   x_{k+1}^{\pm}(s), & 0 \le k \le p-s-2, \\
                   a_0^{\pm}(s), & k=p-s-1,
                   \end{cases}\notag\\
&Fy_{k}^{\pm}(s) = \begin{cases}
                   y_{k+1}^{\pm}(s), & 0 \le k \le p-s-2, \\
                   0, & k=p-s-1,
                   \end{cases}\notag\\  
&Fa_{n}^{\pm}(s) = \begin{cases}
                    a_{n+1}^{\pm}(s), & 0 \le n \le s-2, \\
                    0, & n=s-1.
                    \end{cases}\notag
\end{align}
\item[$(2)$]
The module $\mathcal{P}_s^{\pm}$ is indecomposable.
\end{enumerate}
\end{Prop} 
\begin{proof}
$(1)$We see that
\begin{align}
Kb_n^{\pm}(s) &= \pm q^{s-1-2n}b_n^{\pm}(s),  \label{proj:re1} \\
Ka_n^{\pm}(s) &= \pm q^{s-1-2n}a_n^{\pm}(s),  \label{proj:re2} \\
Kx_{k}^{\pm}(s) &= \mp q^{p-s-1-2k}x_{k}^{\pm}(s), \label{proj:re3} \\
Ky_{k}^{\pm}(s) &= \mp q^{p-s-1-2k} y_{k}^{\pm}(s) \label{proj:re4}
\end{align}
for $\ 0 \le n \le s-1$ and $ \ 0 \le k \le p-s-1$.
By Lemma \ref{proj:3} and (\ref{proj:re0}), we have
\begin{align}
FEb_0^{\pm}(s) = F x_{p-s-1}^{\pm}(s) = a_0^{\pm}(s). \label{proj:re5}
\end{align}
By Lemma \ref{Lem:2}, (\ref{proj:re1}) and (\ref{proj:re5}),
\begin{align}
Ey_{0}^{\pm}(s) &= EF^{s} b_0^{\pm}(s) \notag\\
&= \left\{ F^sE+ \frac{[s]}{q-q^{-1}}F^{s-1}(q^{-(s-1)}K-q^{s-1}K{-1}) \right\} b_0^{\pm}(s) \notag\\
&=  F^sE b_0^{\pm}(s) = a_{s-1}^{\pm}(s) \label{proj:re6}
\end{align}
and this shows
\begin{align}
Ey_{k}^{\pm}(s) &= EF^k y_{0}^{\pm}(s) \notag\\
&=\left\{ F^kE+ \frac{[k]}{q-q^{-1}}F^{k-1}(q^{-(k-1)}K-q^{k-1}K{-1}) \right\} y_{0}^{\pm}(s) \notag\\
&=\mp [k][p-s-k] y_{k-1}^{\pm} \label{proj:re7}
\end{align}
for $1 \le k \le p-s-1$.
Since 
\begin{align}
E^{p-s+1}b_{0}^{\pm}(s) &= E^{p-s+1}\sum_{n=1}^{p-s} \alpha_n^{\pm} (s) E^{p-1-n}F^{p-1-n}v_s^{\pm} \notag\\
&=0, \notag
\end{align}
we have
\begin{align}
Ex_{0}^{\pm}(s) &= \frac{1}{\prod_{\ell = 1}^{p-s-1}(\mp[\ell][p-s-\ell])}E^{p-s}x_{p-s-1}^{\pm}(s) \notag\\
&= \frac{1}{\prod_{i=1}^{s-1}(\pm[i][s-i])\prod_{\ell = 1}^{p-s-1}(\mp[\ell][p-s-\ell])}E^{p-s+1}b_{0}^{\pm}(s) \notag\\
&=0. \label{proj:re11}
\end{align}
Therefore we have
\begin{align}
&F x_{k}^{\pm}(s) \notag\\
&= \frac{1}{\prod_{\ell = k+1}^{p-s-1}(\mp[\ell][p-s-\ell])}FE^{p-s-k}b_{0}^{\pm}(s) \notag\\
&= \frac{1}{\prod_{\ell = k+1}^{p-s-1}(\mp[\ell][p-s-\ell])} \notag\\
&\times \left\{ E^{p-s-1-k}F
 - \frac{[p-s-1-k]}{q-q^{-1}}E^{p-s-2-k}(q^{p-s-2-k}K-q^{-(p-s-2-k)}K^{-1}) \right\} \notag\\
&\times  x_{p-s-1}^{\pm}(s) \notag\\
&= \frac{\mp [k+1][p-s-1-k]}{\prod_{\ell = k+1}^{p-s-1}(\mp[\ell][p-s-\ell])}E^{p-s-1-k}b_{0}^{\pm}(s) \notag\\
&= x_{k-1}^{\pm}(s) \label{proj:re12}
\end{align}
for $1 \le k \le p-s-1$.
The relation (\ref{proj:re6}) shows
\begin{align}
Eb_n^{\pm}(s) &= EF^n b_0^{\pm}(s) \notag\\
&= \left\{ F^nE + \frac{[n]}{q-q^{-1}}F^{n-1}(q^{-(n-1)}K-q^{n-1}K{-1})  \right\} b_0^{\pm}(s) \notag\\
&= \pm [n][s-n] b_{n-1}^{\pm}(s) +   a_{n-1}^{\pm}(s) \label{proj:re8}
\end{align} 
for $1 \le n \le s-1$.
By the definition of the vector $x_{k}^{\pm}(s)$, we obtain
\begin{align}
Ex_{k}^{\pm}(s) &= \mp [k][p-s-k] \frac{1}{\prod_{\ell = k}^{p-s-1}(\mp[\ell][p-s-\ell])}E^{p-s-(k-1)}b_{0}^{\pm}(s) \notag\\
&= \mp [k][p-s-k]x_{k-1}^{\pm}(s) \label{proj:re9}
\end{align}
for $1 \le k \le p-s-1$.
By using (\ref{proj:re8}) inductively, we can show that
\begin{align}
&E^{s-1} b_{s-1}^{\pm}(s) \notag\\
&= \prod_{i=1}^{s-1}(\pm[i][s-i])b_{0}^{\pm}(s) 
+ \sum_{i=1}^{s-1} \prod_{\genfrac{}{}{0pt}{}{j=1}{j \neq i}}^{s-1} (\pm[j][s-j])a_{0}^{\pm}(s). \label{proj:re10}
\end{align}
By the definition of vectors and the relation (\ref{proj:re1})--(\ref{proj:re9}) and (\ref{proj:re11})--(\ref{proj:re12}),
the space $\mathcal{P}_s^{\pm}$ is a $2p$-dimensional left $\overline{U}_q$-module.

$(2)$The module $\mathcal{P}_s^{\pm}$ has unique highest weight vector which generates the simple module $\mathcal{X}_s^{\pm}$
up to scalar multiple,
so the module $\mathcal{P}_s^{\pm}$ is indecomposable.
\end{proof}

Next we show that the generator of $\mathcal{P}_s^{\pm}$ is an idempotent.
Now we can see that
\begin{align}
(b_0^{\pm}(s))^2= \sum_{n=1}^{p-s} \alpha_n^{\pm}(s)E^{p-1-n}F^{p-1-n}v_s^{\pm}b_0^{\pm}(s).\notag
\end{align}
By Proposition \ref{proj:4}, one has
\begin{align}
&\alpha_n^{\pm}(s)E^{p-1-n}F^{p-1-n}v_s^{\pm}b_0^{\pm}(s) \notag\\
&=2p \alpha_n^{\pm}(s) E^{p-1-n} y_{p-s-1-n}^{\pm}(s) \notag\\
&=2p \alpha_n^{\pm}(s) \prod_{i=1}^{s-1} (\pm[i][s-i]) \prod_{k=1}^{p-s-1-n}(\mp[k][p-s-k]) a_0^{\pm}(s) \notag\\
&=2p \prod_{i=1}^{s-1} (\pm[i][s-i]) \prod_{\genfrac{}{}{0pt}{}{k=1}{k \neq n}}^{p-s-1}(\mp[k][p-s-k]) a_0^{\pm}(s) \notag
\end{align}
for $1 \le n \le p-s-1$ and by (\ref{proj:re10})
\begin{align}
&\alpha_{p-s}^{\pm}(s)E^{s-1}F^{s-1}v_s^{\pm}b_0^{\pm}(s) \notag\\
&=2p \alpha_{p-s}^{\pm}(s) E^{s-1} b_{s-1}^{\pm}(s) \notag\\
&=2p \prod_{m=1}^{p-s-1}(\mp[m][p-s-m])  \prod_{i=1}^{s-1}(\pm[i][s-i]) b_0^{\pm}(s) \notag\\
& \quad + 2p \prod_{m=1}^{p-s-1}(\mp[m][p-s-m])\sum_{j=1}^{s-1}\prod_{\genfrac{}{}{0pt}{}{ k =1}{k \neq j}}^{s-1}(\pm[k][s-k]) a_0^{\pm}(s). \notag
\end{align}
Moreover we have
\begin{align}
&b_0^{\pm}(s) a_0^{\pm}(s) \notag\\&
= \alpha_{p-s}^{\pm}(s)E^{s-1}F^{s-1}v_s^{\pm}a_0^{\pm}(s) \notag\\
&= 2p \prod_{m=1}^{p-s-1}(\mp[m][p-s-m]) \prod_{i=1}^{s-1} (\pm[i][s-i]) a_0^{\pm}(s) \notag
\end{align}
and
\begin{align}
&a_0^{\pm}(s) b_0^{\pm}(s) = 2p \prod_{m=1}^{p-s-1}(\mp[m][p-s-m]) \prod_{i=1}^{s-1} (\pm[i][s-i]) a_0^{\pm}(s). \notag
\end{align}
We set
\begin{align}
&\gamma^{\pm}(s) = 2p \prod_{m=1}^{p-s-1}(\mp[m][p-s-m]) \prod_{i=1}^{s-1} (\pm[i][s-i]), \label{eq:gamma}\\
&\delta^{\pm}(s) = 2p \prod_{m=1}^{p-s-1}(\mp[m][p-s-m])\sum_{j=1}^{s-1}\prod_{\genfrac{}{}{0pt}{}{ k =1}{k \neq j}}^{s-1}(\pm[k][s-k])  \notag\\
& \qquad \qquad \qquad + 2p \prod_{i=1}^{s-1} (\pm[i][s-i])\sum_{n=1}^{p-s-1}\prod_{\genfrac{}{}{0pt}{}{k=1}{k \neq n}}^{p-s-1}(\mp[k][p-s-k]) \label{eq:delta}
\end{align}
and
\begin{align}
e_{s}^{\pm} = \frac{1}{\gamma^{\pm}(s)} \left(b_0^{\pm}(s) - \frac{\delta^{\pm}(s)}{\gamma^{\pm}(s)}a_0^{\pm}(s) \right),\notag
\end{align}
then we have
\begin{align}
&(b_0^{\pm}(s))^2 = \gamma^{\pm}(s)b_0^{\pm}(s) + \delta^{\pm}(s)a_0^{\pm}(s), \label{eq:proj1}\\ 
&a_0^{\pm}(s) b_0^{\pm}(s) = \gamma^{\pm}(s) a_0^{\pm}(s), \label{eq:proj2}\\ 
&b_0^{\pm}(s)a_0^{\pm}(s) = \gamma^{\pm}(s) a_0^{\pm}(s). \label{eq:proj3}
\end{align}
By Proposition \ref{proj:4} we can see that the vectors
\begin{align}
&\frac{1}{\gamma^{\pm}(s)} \left(b_n^{\pm}(s) - \frac{\delta^{\pm}(s)}{\gamma^{\pm}(s)}a_n^{\pm}(s) \right), \notag\\
&\frac{1}{\gamma^{\pm}(s)} x_k^{\pm}(s), \notag\\
&\frac{1}{\gamma^{\pm}(s)} x_k^{\pm}(s), \notag\\
&\frac{1}{\gamma^{\pm}(s)} a_n^{\pm}(s), \notag
\end{align}
for $0 \le n \le s-1$ and $0 \le k \le p-s-1$
also form the basis of $\mathcal{P}_s^{\pm}$.

\begin{Prop}\label{proj:6}
The element $e_{s}^{\pm}$ is a primitive idempotent of $\overline{U}_q$.
In paticular, the module $\mathcal{P}_s^{\pm}$ is projective.
\end{Prop}
\begin{proof}
We see $(a_0^{\pm}(s))^2 = 0$, then by $(\ref{eq:proj1})$--$(\ref{eq:proj3})$,
\begin{align}
(e_{s}^{\pm})^{2} &= \frac{1}{(\gamma^{\pm}(s))^2} \{ (b_0^{\pm}(s))^2 - \frac{\delta^{\pm}(s)}{\gamma^{\pm}(s)}a_0^{\pm}(s)b_0^{\pm}(s)
-\frac{\delta^{\pm}(s)}{\gamma^{\pm}(s)}b_0^{\pm}(s)a_0^{\pm} \} \notag\\
&= \frac{1}{(\gamma^{\pm}(s))^2} \{ \gamma^{\pm}(s)b_0^{\pm}(s) + \delta^{\pm}(s)a_0^{\pm}(s) - 2 \delta^{\pm}(s)a_0^{\pm}(s) \} \notag\\
&= e_{s}^{\pm}. \notag
\end{align}
By Proposition \ref{proj:4}, the module $\mathcal{P}_s^{\pm}$ is generated by the idempotent $e_{s}^{\pm}$. 
\end{proof}
By Proposition \ref{proj:2} and Proposition \ref{proj:6}, we obtain the following result:
\begin{Thm}\label{proj:7}
$\{e_s^{\pm} | \ 1 \le s \le p \}$ is a basic set of primitive idempotents of $\overline{U}_q$.
\end{Thm} 

\subsection{The Casimir element}
Let us denote the Casimir element $C$ defined by
\begin{align}
C = EF + \frac{q^{-1}K + qK^{-1}}{(q - q^{-1})^2} = FE + \frac{qK + q^{-1}K^{-1}}{(q - q^{-1})^2}. \notag
\end{align}
Then we can see 
\begin{align}
C e_s^{\pm} &= \left\{ FE + \frac{qK + q^{-1}K^{-1}}{(q - q^{-1})^2} \right\} e_s^{\pm} \notag\\
& = a_0^{\pm} (s) \pm \frac{q^s + q^{-s}}{(q - q^{-1})^2} e_s^{\pm} \notag
\end{align}
for $1 \le s \le p-1$ and
\begin{align}
Ca_0^{\pm} (s) =  \pm \frac{q^s + q^{-s}}{(q - q^{-1})^2} a_0^{\pm}(s) \notag
\end{align}
for $1 \le s \le p$.
Thus we have the following proposition.
\begin{Prop}[\cite{FGST1}]
The minimal polynomial of $C$ is given by 
\begin{align}
\Psi_{2p}(x) = (x - \beta_0 )(x - \beta_p ) \prod_{j=1}^{p-1} (x - \beta_j )^2, \ \beta_j = \frac{q^j+q^{-j}}{(q - q^{-1})^2}. \notag
\end{align} 
\end{Prop}
By this proposition, $\overline{U}_q$ has the decomposition into subalgebras
\begin{align}
\overline{U}_q = \bigoplus_{s=0}^{p} \mathcal{Q}_s \label{eq:subalgebra}
\end{align}
where $\mathcal{Q}_s$ has the eigenvalue $\beta_s$ of $C$.
Moreover we have $e_s^{+}, e_{p-s}^{-} \in \mathcal{Q}_s$, $e_p^- \in \mathcal{Q}_0$ and $e_p^+ \in \mathcal{Q}_p$.

\section{The space of symmetric linear functions of $\overline{U}_q$}
In this section we determine the dimension of the space of symmetric linear functions of $\overline{U}_q$.

\subsection{The basic algebra of $\overline{U}_q$}
Recall that the basic set of  primitive idempotents of $\overline{U}_q$ is given by
\begin{align}
&e_s^{\pm}= \frac{1}{\gamma^{\pm} (s)}(\sum_{n=1}^{p-s} \alpha_n^{\pm} (s) E^{p-1-n}F^{p-1-n}- \frac{\delta^{\pm} (s)}{\gamma^{\pm} (s)} E^{p-1}F^{p-1})v_s^{\pm}, \notag\\
&e_p^{\pm}= \frac{1}{2p \prod_{i=1}^{p-1}(\pm[i][s-i])}a_0^{\pm}(p) \notag
\end{align}
for $1 \le s \le p-1$.
By $(\ref{eq:gamma})$ and $(\ref{eq:delta})$, we have
\begin{align}
\gamma = \gamma^+(s) = 2p \prod_{m=1}^{p-s-1}(-[m][p-s-m]) \prod_{i=1}^{s-1} [i][s-i]=\gamma^-(p-s), \notag
\end{align}
and
\begin{align}
\delta 
=\delta^{+}(s) &= 2p \prod_{m=1}^{p-s-1}(-[m][p-s-m])\sum_{j=1}^{s-1}\prod_{\genfrac{}{}{0pt}{}{ k =1}{k \neq j}}^{s-1}[k][s-k]  \notag\\
& \qquad \qquad \qquad + 2p \prod_{i=1}^{s-1} [i][s-i]\sum_{n=1}^{p-s-1}\prod_{\genfrac{}{}{0pt}{}{k=1}{k \neq n}}^{p-s-1}(-[k][p-s-k]), \notag\\
&=\delta^{-}(p-s).\notag
\end{align}
Thus we have
\begin{align}
&e_s^{+}= \frac{1}{\gamma}(\sum_{n=1}^{p-s} \alpha_n^{+} (s) E^{p-1-n}F^{p-1-n}- \frac{\delta}{\gamma} E^{p-1}F^{p-1})v_s^{+}, \notag\\
&e_{p-s}^- = \frac{1}{\gamma}(\sum_{n=1}^{s} \alpha_n^{-} (p-s) E^{p-1-n}F^{p-1-n}- \frac{\delta}{\gamma} E^{p-1}F^{p-1})v_{p-s}^{-}. \notag
\end{align}

The basic algebra of $\overline{U}_q$ is defined by $e \overline{U}_q e$ where 
$e= \sum_{s=1}^{p-1} (e_s^+ + e_{p-s}^-) + e_p^+ + e_p^-$.
By the decomposition of $\overline{U}_q$, we can see
\begin{align}
e\overline{U}_qe = \bigoplus_{s=1}^{p-1}(e_s^+ + e_{p-s}^-)\mathcal{Q}_s(e_s^+ + e_{p-s}^-) \oplus e_p^+ \mathcal{Q}_p e_p^+ \oplus e_p^- \mathcal{Q}_0 e_p^-. \notag
\end{align}
We can see that $(e_s^+ + e_{p-s}^-)\mathcal{Q}_s(e_s^+ + e_{p-s}^-)$ is the basic algebra of $\mathcal{Q}_s$ and each of $e_p^+ \mathcal{Q}_p e_p^+$ and $e_p^- \mathcal{Q}_0 e_p^-$
is the basic algebra of $\mathcal{Q}_p$ and $\mathcal{Q}_0$, respectively.
We set $B_s = (e_s^+ + e_{p-s}^-)\mathcal{Q}_s(e_s^+ + e_{p-s}^-)$ for $1 \le s \le p-1$, $B_p = e_p^+ \mathcal{Q}_p e_p^+$ and $B_0 = e_p^- \mathcal{Q}_0 e_p^-$.
By the definition of $e_p^{\pm}$ and the structure of the simple module $\mathcal{X}_p^{\pm}$,
we have the following: 
\begin{Prop}\label{Prop:cobasic}
The basic algebras $B_0$ and $B_p$ are both $1$-dimensional commutative algebras.
\end{Prop}

In order to determine the structure of the basic algebra $B_s$ for $1 \le s \le p-1$,
we make the table of multiplications of $B_s$.

We can see that $\mathcal{Q}_s(e_s^+ + e_{p-s}^-) \cong \mathcal{P}_s^+ \oplus \mathcal{P}_{p-s}^-$.
Then the basis of $\mathcal{Q}_s(e_s^+ + e_{p-s}^-)$ corresponding to $\mathcal{P}_s^+$ is given by
\begin{align}
&B_n^{+} = F^n e_s^+, \notag\\
&Y_k^{+} = F^{s+k} e_s^+ = \frac{1}{\gamma} y_k^+(s), \notag\\
&X_k^{+} = \frac{1}{\prod_{\ell = k+1}^{p-s-1}(-[\ell][p-s-\ell])}E^{p-s-k} e_s^+ = \frac{1}{\gamma}x_k^+(s),\notag\\
&A_n^{+} = F^{n+1}E e_s^+ = \frac{1}{\gamma} a_n^{+}(s) \notag
\end{align} 
and that corresopnding to $\mathcal{P}_{p-s}^-$ is given by
\begin{align}
&B_k^{-} = F^k e_{p-s}^-, \notag\\
&Y_n^{-} = F^{p-s+n} e_{p-s}^- = \frac{1}{\gamma} y_n^-(p-s), \notag\\
&X_n^{-} = \frac{1}{\prod_{\ell = n+1}^{s-1}[\ell][s-\ell]}E^{s-k}  e_{p-s}^- = \frac{1}{\gamma}x_n^-(p-s), \notag\\
&A_k^{-} = F^{k+1}E e_{p-s}^- = \frac{1}{\gamma} a_n^{-}(p-s) \notag
\end{align}
for $0 \le n \le s-1$ and $0 \le k \le p-s-1$.

The basis of $B_s$ is the set of non-vanishing elements of the basis of $\mathcal{Q}_s(e_s^+ + e_{p-s}^-)$
by the left action of primitive idempotents $e_s^+$ and $e_{p-s}^-$.
For finding such the elements, we need the follwing lemma:

\begin{Lem}\label{Lem:eigen}
Let $\varphi_n$ be a vector of weight $q^{s-1-2n}$ and $\psi _k$
be a vector of weight $-q^{p-s-1-2k}$ for $0 \le n \le s-1$ and $0 \le k \le p-s-1$.
Then
\begin{align}
v_s^+ \varphi_n &= \begin{cases}
                  2p\varphi _0, & n=0, \\
                  0, & n \neq 0, 
                 \end{cases}\notag\\
v_s^+  \psi _k&= 0, \notag\\
v_{p-s}^- \varphi_n &= 0, \notag\\
v_{p-s}^- \psi _k&= \begin{cases}
                  2p\psi _0, & k=0, \\
                  0, & k \neq 0. 
                 \end{cases}\notag
\end{align}
\end{Lem}
\begin{proof}
Since $v_s^+ = \sum_{\ell = 0}^{2p-1} q^{-(s-1) \ell} K^{\ell}$, we have 
\begin{align}
v_s^+ \varphi_n &= \sum_{\ell = 0}^{2p-1} q^{-(s-1) \ell} (q^{s-1-2n})^{\ell}\varphi_n \notag\\
& = \sum_{\ell = 0}^{2p-1} q^{-2n \ell} \varphi_n. \notag
\end{align}
We can see that $1-q^{-2n} = 0$ if and only if $n=0$ since $0 \le n \le s-1$.
\end{proof}
By Lemma \ref{Lem:eigen}, we can get
\begin{align}
&e_s^+ B_n^{+} =0, \ e_s^+ Y_k^{+} =0, \ e_s^+ X_{k}^{+} =0, \ e_s^+ A_{n}^{+} =0, \notag\\
&e_s^+ B_k^{-} =0, \ e_s^+ Y_n^{-} =0, \ e_s^+ X_{n}^{-} =0, \ e_s^+ A_{k}^{-} =0, \notag
\end{align}
for $1 \le n \le s-1$ and $0 \le k \le p-s-1$ and
\begin{align}
&e_{p-s}^- B_m^{+} =0, \ e_{p-s}^- Y_j^{+} =0, \ e_{p-s}^- X_{j}^{+} =0, \ e_{p-s}^- A_{m}^{+} =0, \notag\\
&e_{p-s}^- B_j^{-} =0, \ e_{p-s}^- Y_m^{-} =0, \ e_{p-s}^- X_{m}^{-} =0, \ e_{p-s}^- A_{j}^{-} =0, \notag
\end{align}
for $0 \le m \le s-1$ and $1 \le j \le p-s-1$.
By Lemma \ref{Lem:eigen} and the definition of idempotents, we get 
\begin{align}
e_s^+ Y_0^{-} &= \frac{2p}{\gamma}(\sum_{n=1}^{p-s} \alpha_n^+ (s) E^{p-1-n}F^{p-1-n}- \frac{\delta }{\gamma } E^{p-1}F^{p-1})Y_0^{-} \notag\\
            &= \frac{2p\alpha^+_{p-s} (s)}{\gamma}E^{s-1}F^{s-1}Y_0^- \notag\\
            &= Y_0^-, \label{basic:re1} \\
e_s^+ X_0^{-} &= \frac{2p}{\gamma}(\sum_{n=1}^{p-s} \alpha_n^+ (s) E^{p-1-n}F^{p-1-n}- \frac{\delta}{\gamma} E^{p-1}F^{p-1})X_0^{-} \notag\\
            &= \frac{2p\alpha^+_{p-s} (s)}{\gamma}E^{s-1}F^{s-1}X_0^- \notag\\
            &= X_0^-, \label{basic:re2} \\
e_s^+ A_0^{+} &= \frac{2p}{\gamma}(\sum_{n=1}^{p-s} \alpha_n^+ (s) E^{p-1-n}F^{p-1-n}- \frac{\delta}{\gamma} E^{p-1}F^{p-1})A_0^{+} \notag\\
            &= \frac{2p\alpha_{p-s}^+ (s)}{\gamma}E^{s-1}F^{s-1}A_0^+ \notag\\
            &= A_0^+. \label{basic:re3}
\end{align}
Similarly we have
\begin{align}
e_{p-s}^- X_0^+ &= X_0^+, \label{basic:re4} \\
e_{p-s}^- Y_0^+ &=Y_0^+, \label{basic:re5} \\
e_{p-s}^- A_0^- &= A_0^-. \label{basic:re6}
\end{align}
This implies that the basic algebra $B_s$ for $1 \le s \le p-1$ has the basis of the form
\begin{align}
e_s^{+}, \ e_{p-s}^-, \ X_0^{\pm}, \ Y_0^{\pm}, \ A_0^{\pm}. \notag
\end{align}

Next we determine the multiplication table of $B_s$.
By Lemma \ref{Lem:eigen} and the definition of the vectors, we can see that
for $a = X_0^+, Y_0^+$, or $A_0^{+}$,
\begin{align}
a b = \begin{cases}
                    a, & b= e_s^+, \\
                    0, & b = X_0^+, Y_0^+, e_{p-s}^-,  A_0^-, 
                    \end{cases}\notag
\end{align}
and for $a = X_0^-, Y_0^-$, or $A_0^{-}$,
\begin{align}
a b = \begin{cases}
                    a, & b= e_{p-s}^-, \\
                    0, & b = X_0^-, Y_0^-, e_s^+,  A_0^+.
                    \end{cases}\notag
\end{align}
By (\ref{basic:re1}), we have
\begin{align}
X_0^+Y_0^- &= \frac{1}{\prod_{\ell = 1}^{p-s-1}(-[\ell][p-s-\ell])}E^{p-s} e_s^+Y_0^{-} \notag\\
           &= \frac{1}{\prod_{\ell = 1}^{p-s-1}(-[\ell][p-s-\ell])}E^{p-s} Y_0^{-} \notag\\
           &= A_0^-, \notag\\
Y_0^+Y_0^- &= F^{s} e_s^+ Y_0^- = 0, \notag\\
A_0^+Y_0^- &= FEe_s^+ Y_0^- = FA_{p-s-1}^- = 0, \notag
\end{align}
and by (\ref{basic:re2}),
\begin{align}
X_0^+X_0^- &= \frac{1}{\prod_{\ell = 1}^{p-s-1}(-[\ell][p-s-\ell])}E^{p-s} e_s^+X_0^- \notag\\
           &= 0 \notag\\
Y_0^+X_0^- &= F^{s} e_s^+ X_0^- = A_0^-, \notag\\
A_0^+X_0^- &= FEe_s^+ X_0^- = 0, \notag
\end{align}
and by (\ref{basic:re3}),         
\begin{align}
X_0^+A_0^+ &= \frac{1}{\prod_{\ell = 1}^{p-s-1}(-[\ell][p-s-\ell])}E^{p-s} e_s^+A_0^+ \notag\\
           &= 0, \notag\\
Y_0^+A_0^+ &= F^{s} e_s^+ A_0^+ = 0, \notag\\
A_0^+A_0^+ &= FEe_s^+ A_0^+ = 0. \notag
\end{align}
By (\ref{basic:re4})--(\ref{basic:re6}) and similar caluculation appeared above, we get
\begin{align}
a  X_0^+&= \begin{cases}
                    A_0^+, & a=Y_0^{-}, \\
                    0, & a = X_0^-,A_0^-,
                    \end{cases}\notag\\
a  Y_0^+&= \begin{cases}
                    A_0^+, & a=X_0^{-}, \\
                    0, & a = Y_0^-,A_0^-,
                    \end{cases}\notag
\end{align}
and
\begin{align}
a  A_0^-= 0. \notag
\end{align}
for $a =  X_0^-, Y_0^{-}, A_0^-$.

The above relations yield the following multipulication table among the basis of $B_s$ for $1 \le s \le p-1$:
\begin{center}
Table 1 : The table of $xy$ for $x, y = e_s, e_{p-s}^-, A_0^{\pm}, X_0^{\pm}, Y_0^{\pm}$ 

\begin{tabular}{c|cccccccc}
\multicolumn{9}{c}{} \\
$x \backslash y$    & $e_s^+$   & $X_0^+$ & $Y_0^+$ & $A_0^+$ & $e_{p-s}^-$   & $X_0^-$ & $Y_0^-$ & $A_0^-$ \\ \hline
$e_s^+$   & $e_s^+$   & $0$     & $0$     & $A_0^+$ & $0$     & $X_0^-$ & $Y_0^-$ & $0$     \\ 
$X_0^+$ & $X_0^+$ & $0$     & $0$     & $0$     & $0$     & $0$     & $A_0^-$ & $0$     \\ 
$Y_0^+$ & $Y_0^+$ & $0$     & $0$     & $0$     & $0$     & $A_0^-$ & $0$     & $0$     \\ 
$A_0^+$ & $A_0^+$ & $0$     & $0$     & $0$     & $0$     & $0$     & $0$     & $0$     \\ 
$e_{p-s}^-$   & $0$     & $X_0^+$ & $Y_0^+$ & $0$     & $e_{p-s}^-$   & $0$     & $0$     & $A_0^-$ \\ 
$X_0^-$ & $0$     & $0$     & $A_0^+$ & $0$     & $X_0^-$ & $0$     & $0$     & $0$     \\ 
$Y_0^-$ & $0$     & $A_0^+$ & $0$     & $0$     & $Y_0^-$ & $0$     & $0$     & $0$     \\ 
$A_0^-$ & $0$     & $0$     & $0$     & $0$     & $A_0^-$ & $0$     & $0$     & $0$     \\ 
\end{tabular}
\end{center}
\begin{Thm}\label{Thm:basic1}
The basic algebra $B_s$ for $1 \le s \le p-1$ is $8$-dimensional associative algebra with the basis
\begin{align}
e_s^{+}, \ e_{p-s}^-, \ X_0^{\pm}, \ Y_0^{\pm}, \ A_0^{\pm}, \notag
\end{align}
and the relations in Table $1$.
\end{Thm}

\subsection{The dimension of the space of symmetric linear functions of $\overline{U}_q$}
In this section, we determine the dimension
of the space of symmetric linear functions of $\overline{U}_q$.

First we determine the subspace $[B_s,B_s]$ of $B_s$ by
Proposition \ref{Prop:cobasic} and Theorem \ref{Thm:basic1}.
By Table 1, we can obtain all the commutators among the basis of $B_s$ for $1 \le s \le p-1 $,
in particular, we can see $[x, A_0^{\pm}] = 0$ for $x = e_s^+, e_{p-s}^-, X_0^{\pm}, Y_0^{\pm}, A_0^{\pm}$.
\begin{center}
\samepage
Table 2 : The table of $[x, y]$ for $x, y = e_s^+, e_{p-s}^-, X_0^{\pm}, Y_0^{\pm}$ 
\begin{tabular}{c|cccccc}
\multicolumn{7}{c}{} \\
$x \backslash y$    & $e_s^+$ & $X_0^+$       & $Y_0^+$         & $e_{p-s}^-$ & $X_0^-$         & $Y_0^-$   \\ \hline
$e_s^+$             & $0$     & $-X_0^+$      & $-Y_0^+$        & $0$         & $X_0^-$         & $Y_0^-$       \\ 
$X_0^+$             & $X_0^+$ & $0$           & $0$             & $-X_0^+$    & $0$             & $A_0^- - A_0^+$      \\ 
$Y_0^+$             & $Y_0^+$ & $0$           & $0$             & $-Y_0^+$    & $A_0^- - A_0^+$ & $0$         \\ 
$e_{p-s}^-$         & $0$     & $X_0^+$       & $Y_0^+$         & $0$         & $-X_0^-$        & $-Y_0^-$     \\ 
$X_0^-$             & $-X_0^-$& $-X_0^-$      & $A_0^+-A_0^-$   & $X_0^-$     & $0$             & $0$          \\ 
$Y_0^-$             & $-Y_0^-$& $A_0^+-A_0^-$ & $0$             & $Y_0^-$     & $0$             & $0$\\ 
\end{tabular}
\end{center}

\begin{Prop}\label{Prop:basicsym}
The space of symmetric linear functions of $B_s$ is $3$-dimensional for $1 \le s \le p-1 $
and $1$-dimensional for $s = 0, p$.
\end{Prop}
\begin{proof}
The space $B_0/[B_0, B_0]$ and $B_p/[B_p, B_p]$ are $1$-dimensional
since $B_0$ and $B_p$ are $1$-dimensional commutative algebras.

Table 2 shows that $\{ X_0^{\pm}, Y_0^{\pm}, A_0^+-A_0^- \}$ is the basis of the space $[B_s, B_s]$.
Thus the dimension of the space $B_s/[B_s, B_s]$ is equal to $3$ for $1 \le s \le p-1 $. 
\end{proof}

By this proposition, we have the following:
\begin{Thm}\label{Thm:main}
The space of symmetric linear functions of the basic algebra of $\overline{U}_q$ is $(3p-1)$-dimensional.
\end{Thm}
By Theorem \ref{Thm:NS} and Theorem \ref{Thm:main}, we finally get the following theorem:
\begin{Thm}
The space of symmetric linear functions of $\overline{U}_q$ is $(3p-1)$-dimensional.
\end{Thm}

\end{document}